----------
X-Sun-Data-Type: tex-file
X-Sun-Data-Description: tex-file
X-Sun-Data-Name: toroidal.tex
X-Sun-Charset: iso-8859-1
X-Sun-Content-Lines: 1236

\documentclass[11pt]{article}
\usepackage{epsfig}
\usepackage{amsmath}
\usepackage{amsfonts}
\usepackage{amssymb}

\numberwithin{equation}{section}

\title{Twisted Toroidal Lie Algebras}
\author{ Johan van de Leur\\
\\
Mathematical Institute,\\
University of Utrecht,\\
P.O. Box 80010, 3508 TA Utrecht,\\
The Netherlands\\
e-mail: vdleur@math.uu.nl}

\newtheorem{definition}{Definition}[section]
\newtheorem{lemma}{Lemma}[section]

\begin{document}

\maketitle

\begin{abstract}
Using $n$ finite order automorphisms on a simple complex Lie algebra 
we construct twisted $n$-toroidal Lie algebras. Thus obtaining Lie algebras wich have a rootspace decomposition.
For the case $n=2$ we list certain  simple Lie algebras 
and their automorphisms, which produce twisted 2-toroidal algebras.  In this way we obtain Lie algebras that are related to all Extended 
Affine Root Systems of K. Saito.
\end{abstract}

\section{Introduction}

At the end of the 1960's Victor Kac \cite{K1} and Bob Moody \cite{M}
independently realized that one could generalize Serre's construction of simple
Lie algebras to  construct certain infinite-dimensional Lie algebras, that
posess root systems. A special class of these Kac-Moody (Lie) algebras were
related to the affine root systems. These affine Lie algebras could also be
constructed in an explicit way. Nowadays this construction is well known. These
algebras are the central extensions of so-called twisted and untwisted  loop
algebras. The untwisted loop algebras
$\tilde\mathfrak{g}$ can be obtained as follows.
Let $\dot\mathfrak{g}$ be a finite dimensional simple Lie algebra, then
$\tilde\mathfrak{g}=\dot\mathfrak{g}\otimes\mathbb{C}[t,t^{-1}]$, where
$\mathbb{C}[t,t^{-1}]$ is the algebra of Laurent polynomials in the variable
$t$. The twisted loop algebras can be obtained as certain subalgebras of the
untwisted ones. Let $\sigma$ be a certain finite order automorphism induced by
a diagram automorphism of the Dynkin diagram (see \cite{Kac} for more details).
Then $\dot\mathfrak{g}$ decomposes into eigenspaces with respect to this
automorphism. To be more precise, let $n$ be the order of the automorphism and
$\epsilon= e^{\frac {2\pi i}{n}}$
\[
\begin{aligned}
\dot{\mathfrak{g}} &=\bigoplus_{\overline k\in \mathbb{Z}/n\mathbb{Z}}
\dot{\mathfrak{g}}_{\overline k}, \quad\text{where}\\[2mm]
\dot{\mathfrak{g}}_{\overline k}&=\{ g\in
\dot{\mathfrak{g}}|\sigma(g)=\epsilon^{k} g \}.
\end{aligned}
\]
The twisted loop algebra is the following subalgebra of $\tilde\mathfrak{g}$
\[
{\tilde\mathfrak{g}}(\sigma)=\bigoplus_{k\in\mathbb{Z}}
\dot{\mathfrak{g}}_{\overline k}\otimes t^{k}.
\]
A generalization of this construction, at least of the untwisted ones is
clearly obvious. Instead of tensoring by the algebra of Laurent polynomials in
one variable, one can take Laurent polynomials in $N$ variables $t_j$. Thus
obtaining toroidal Lie algebras. Unfortunately, these Lie algebras are not
Kac-Moody algebras, but they still are very interesting and obviously related
to certain extensions of affine root systems. They appeared in the work of
Slodowy \cite{Sl2} as certain intersection matrix algebras.

K. Saito, interested in singularity theory and inspired by the work of
Looijenga \cite{L1}, \cite{L2} and Slodowy \cite{Sl1}, \cite{Sl2}, classified
in \cite{S1} extended affine root systems, whose radical is 2-dimensional and
for which the quotient of the root system modulo a certain 1-dimensional space
is reduced.
In 1997 Allison, Azam, Berman, Gao and Pianzola \cite{AA} had a different
approach, they  classified these root systems using semilattices.
The Lie algebras corresponding to these root systems, so-called extended affine
Lie algebras or more precise bi-affine Lie algebras, were constructed in a
paper by Hoegh--Krohn and Torresani \cite{HT}. However, their construction was
not complete. Pollmann gave a complete construction in \cite{P}, which was
based on the idea of twisting affine Lie algebras by finite order automorphism.
This idea was presented in \cite{HT} and also in an unpublished paper of
Wakimoto \cite{W}. Although one obtains in this way the bi-affine algebras as
subalgebras of 2-toroidal Lie algebras, the grading with respect to the
variables of the Laurent polynomials is not so nice. The present paper
constructs the same Lie algebras also as subalgebras of 2-toroidal Lie
algebras, but in a slightly different way. Whereas \cite{HT}, \cite{P} and
\cite{W}
use finite order automorphisms of affine Lie algebras to construct the
bi-affine Lie algebras,  we use two finite order automorphisms of a simple
finite dimensional Lie algebra, which  commute and thus are simultaniously
diagonalizable, to construct them. This construction generalizes in a different
way than the construction of \cite{P} the construction of the twisted affine
Lie algebras. This is an experimental fact, unfortunately, at this moment there
is no general theory or classification which explains this phenomenon. Note
that in some cases one uses only inner automorphisms. Our construction also
easily generalizes to twisted $N$-toroidal Lie algebras, whereas the
generalization of the method of \cite{HT}, \cite{P} and \cite{W} is somewhat
more complicated. A general, but different and  more abstract,   construction
of these extended affine Lie algebras  is given in the AMS Memoir  \cite{AA}
and in \cite{Az}.

The theory of vertex operator constructions of untwisted toroidal Lie algebras
is well developed \cite{BBS}, \cite{Bi}, \cite{RM}, \cite{MEY}, \cite{SY},
\cite{Y}  and applied to hierarchies of soliton equations \cite{Bi1},
\cite{IT}, \cite{ISW1}, \cite{ISW2}. We hope that the construction given in
Section \ref{TA} can be used to define vertex operator constructions on twisted
bi-affine and extended affine algebras.

Untwisted toroidal Lie algebras appear as current algebras of the symmetry of
K\"ahler--Wess--Zumino--Witten models \cite{IKUX}, \cite{IKU}. This is an
extension of (2-dimensional) Wess--Zumino--Witten models on a $2n$-dimensional
K\"ahler manifold. As such it is one possible candidate of a construction of
integrable quantum field theories in more than two dimensions.

The decomposition of some of the exceptional Lie algebras with respect to the
automorphisms  were checked by Willem de Graaf
using \emph{ {GAP -- Groups, Algorithms, and Programming}} \cite{GAP4}.   It is
 a pleasure to thank him and  Prof. P. Slodowy. The latter  for sending the
manuscript \cite{P}.

\section{Extended Affine Root Systems}
\label{EARS}

The following definition of an extended reduced root system can be found in
\cite{S1}, this definition is different from the one in \cite{AA}, there also
the isotropic roots and $0$ are included in the definition.
\begin{definition}
\label{D1}
Let $V$ be a finite dimensional real vector space with a positive semidefinite
symmetric bilinear form $(\cdot,\cdot)_V$.
A subset $R$  of $V$, with
is called an {\rm extended reduced root system} in $V$ if $R$ satisfies the
following axioms:
\begin{itemize}
\item The additve subgroup
$Q(R)=\sum_{\alpha\in R}\mathbb{Z}\alpha$ of $V$ is a full lattice of $V$, i.e.
$R\otimes_\mathbb{R}Q(R)\simeq V$,
\item $(\alpha,\alpha)_V\ne 0$ for any $\alpha\in R$,
\item If $\alpha\in R$, then $2\alpha\not\in R$,
\item For any $\alpha\in R$, $W_\alpha(R)=R$, where
\[
w_\alpha(\beta)=\beta-2\frac{(\beta,\alpha)_V}{(\alpha,\alpha)_V}\alpha,
\]
\item If $\alpha, \beta\in R$, then
\[
 2\frac{(\beta,\alpha)_V}{(\alpha,\alpha)_V}\in \mathbb{Z}
\]
\item $R$ cannot be decomposed in a disjoint union $R_1\cap R_2$, where
$R_1,\ R_2\subset R$, both nonempty, satisfying $(R_1,R_2)_V=0$,
\end{itemize}
\end{definition}
The dimension $\nu$ of the radical
\[
V^0=\{v\in V| (v,w)_V=0 \ \text{for all } w\in V\},
\]
is called the {\it nullity} of the root system $R$.
For $\nu=2$, Saito \cite{S1} classified all marked extended affine root
systems.
Roughly speaking he considered a 1-dimensional marking, which is a linear
subspace $W\subset V^0$ and considered the induced space $V/W$ and
corresponding induced root system. Now assuming that this induced root system
is a reduced (possibly affine) root system, he obtained the following list. See
\cite{S1} or \cite{P} for a more precise statement.
We write $R(X_\ell)$ for the root system of a finite  type $X_\ell$.
\begin{enumerate}
\item $X_\ell^{(1,1)}$, where $X_\ell$ is of type $A\ell,\ B_\ell,\ C_\ell,\
D_\ell,\ E_\ell,\ F_4$ or $G_2$:
\[
R=\{\alpha+m\delta_0+n\delta_1|\alpha\in R(X_\ell),\ n,m\in\mathbb{Z}\}.
\]
\item $X_\ell^{(1,t)}$, where $t=2$ for $X_\ell=B_\ell,\ C_\ell$ and $F_4$, and
 $t=3$ for $X_\ell=G_2$:
\[
\begin{aligned}
R=&\{\alpha+m\delta_0+n\delta_1|\alpha\in R(X_\ell)\ \text{short} ,\
n,m\in\mathbb{Z}\}\\
\ &\cup\{\alpha+m\delta_0+tn\delta_1|\alpha\in R(X_\ell)\ \text{long} ,\
n,m\in\mathbb{Z}\}.
\end{aligned}
\]
\item $X_\ell^{(t,t)}$,where $t=2$ for $X_\ell=B_\ell,\ C_\ell$ and $F_4$, and
$t=3$ for $X_\ell=G_2$:
\[
\begin{aligned}
R=&\{\alpha+m\delta_0+n\delta_1|\alpha\in R(X_\ell)\ \text{short} ,\
n,m\in\mathbb{Z}\}\\
\ &\cup\{\alpha+tm\delta_0+tn\delta_1|\alpha\in R(X_\ell)\ \text{long} ,\
n,m\in\mathbb{Z}\}.
\end{aligned}
\]
\item $A_1^{(1,1)*}$:
\[
R=\{\alpha+m\delta_0+n\delta_1|\alpha\in R(A_1),\ n,m\in\mathbb{Z},\ nm\in
2\mathbb{Z}\}.
\]
\item $B_\ell^{(2,2)*}$:
\[
\begin{aligned}
R=&\{\alpha+m\delta_0+n\delta_1|\alpha\in R(B_\ell)\ \text{short},\
n,m\in\mathbb{Z},\ nm\in 2\mathbb{Z}\}\\
\ &\cup\{\alpha+2m\delta_0+2n\delta_1|\alpha\in R(B_\ell)\ \text{long} ,\
n,m\in\mathbb{Z}\}.
\end{aligned}
\]
\item $C_\ell^{(1,1)*}$:
\[
\begin{aligned}
R=&\{\alpha+m\delta_0+n\delta_1|\alpha\in R(C_\ell)\ \text{long},\
n,m\in\mathbb{Z},\ nm\in 2\mathbb{Z}\}\\
\ &\cup\{\alpha+m\delta_0+n\delta_1|\alpha\in R(C_\ell)\ \text{short} ,\
n,m\in\mathbb{Z}\}.
\end{aligned}
\]
\item $BC_\ell^{(2,1)}$:
\[
\begin{aligned}
R=&\{\alpha+m\delta_0+n\delta_1|\alpha\in R(B_\ell),\  n,m\in\mathbb{Z}\}\\
\ &\cup\{\alpha+(2m+1)\delta_0+n\delta_1|\alpha\in R(C_\ell)\ \text{long} ,\
n,m\in\mathbb{Z}\}.
\end{aligned}
\]
\item $BC_\ell^{(2,2)}(1)$:
\[
\begin{aligned}
R=&\{\alpha+m\delta_0+n\delta_1|\alpha\in R(B_\ell),\  n,m\in\mathbb{Z}\}\\
\ &\cup\{\alpha+(2m+1)\delta_0+2n\delta_1|\alpha\in R(C_\ell)\ \text{long} ,\
n,m\in\mathbb{Z}\}.
\end{aligned}
\]
\item $BC_\ell^{(2,2)}(2)$:
\[
\begin{aligned}
R=&\{\alpha+m\delta_0+n\delta_1|\alpha\in R(B_\ell)\ \text{short},\
n,m\in\mathbb{Z}\}\\
\ &\cup\{\alpha+m\delta_0+2n\delta_1|\alpha\in R(B_\ell)\ \text{long},\
n,m\in\mathbb{Z}\}\\
\ &\cup\{\alpha+(2m+1)\delta_0+2n\delta_1|\alpha\in R(C_\ell)\ \text{long} ,\
n,m\in\mathbb{Z}\}.
\end{aligned}
\]
\item $BC_\ell^{(2,4)}$:
\[
\begin{aligned}
R=&\{\alpha+m\delta_0+n\delta_1|\alpha\in R(B_\ell)\ \text{short},\
n,m\in\mathbb{Z}\}\\
\ &\cup\{\alpha+m\delta_0+2n\delta_1|\alpha\in R(B_\ell)\ \text{long},\
n,m\in\mathbb{Z}\}\\
\ &\cup\{\alpha+(2m+1)\delta_0+4n\delta_1|\alpha\in R(C_\ell)\ \text{long} ,\
n,m\in\mathbb{Z}\}.
\end{aligned}
\]
\item $X_\ell^{(t,1)}$,where $t=2$ for $X_\ell=B_\ell,\ C_\ell$ and $F_4$, and
$t=3$ for $X_\ell=G_2$:
\[
\begin{aligned}
R=&\{ \alpha+m\delta_0+n\delta_1|\alpha\in R(X_\ell)\ \text{short} ,\
n,m\in\mathbb{Z}\}\\
\ &\cup\{\alpha+tm\delta_0+n\delta_1|\alpha\in R(X_\ell)\ \text{long} ,\
n,m\in\mathbb{Z}\}.
\end{aligned}
\]
\end{enumerate}
If we forget the markings, the root systems of type $X_\ell^{(1,t)}$ are
isomorphic to the systems of type $X_\ell^{(t,1)}$.

\section{Toroidal Algebras}
\label{TA}

Lie algebras corresponding to the extended affine root systems of type
$X_\ell^{(1,1)}$ can be easily constructed as follows. Let $\dot{\mathfrak{g}}$
be a simple finite-dimensional complex Lie
algebra with $( \cdot , \cdot )$ the symmetric non-degenerate invariant Killing
form. Let $\dot R$ be its root system. Choose an integer $N \ge 1$ and
consider the tensor product ${\tilde\mathfrak{g}}=\dot{\mathfrak{g}}\otimes
{\cal R}$ of $\dot{\mathfrak{g}}$ with the algebra of Laurent polynomials in
$N+1$ variables:
\[
{\cal R}=\mathbb{C} [t_0^{\pm 1},
t_1^{\pm 1},\ldots, t_n^{\pm 1}]
\]
The toroidal Lie algebra corresponding to $\dot{\mathfrak{g}}$ is
the universal central extension of $\tilde\mathfrak{g}$.
The explicit construction of this extension, which we will present now,  is
known from the papers
\cite{Kas}, \cite{MEY},
see also \cite {BB}.
Let ${\cal K}= \Omega_{\cal R} ^1/d {\cal R}$ be the space of 1-forms modulo
the exact forms. We write $fdg$ for the element of $\cal K$ corresponding to
the pair of elements $f,g $ from $\cal R$ and denote  $k_i =t_i^{-1}d t_i$.
Thus, $\cal K$ is spanned by elements of the form
\[
{t^{\bf m}k_i}=t_0^{m_0}t_1^{m_1} \dots t_N^{m_N}k_i,\quad 0 \leq i \leq N,
\]
 where ${\bf m}= (m_0, m_1, ...,m_N) \in \mathbb{Z}^{N+1}$. Exactness implies
that these elements are related by
\begin{equation}\label{exact}
\sum_{p=0}^{ N} m_p t^{\bf m}k_p =0 , \quad {\bf m} \in \mathbb{Z}^{N+1}.
\end{equation}
Then the toroidal Lie algebra is the vector space
$\hat\mathfrak{g}=\tilde\mathfrak{g}\oplus \cal K$ whith Lie bracket:
\begin{equation}\label{bracket}
[g_1 \otimes f_1(t), g_2\otimes f_2(t)] = [g_1, g_2] \otimes f_1(t) f_2(t) +
(g_1 , g_2) {f_2 d(f_1)} .
\end{equation}
Now, if $N=1$, this clearly gives a Lie algebra whose root system is of type
$X_\ell^{(1,1)}$, viz., the root space corresponding to
$m\delta_0+n\delta_1+\alpha$, with $\alpha\in\dot R$
 is $\hat\mathfrak{g}_{m\delta_0+n\delta_1+\alpha}
=\dot\mathfrak{g}_\alpha\otimes t_0^mt_1^n$.

It is sometimes usefull to add certain outer derivations to the algebra
$\hat\mathfrak{g}$. To do that, we consider
the following algebra of derivations:
\begin{equation}\label{D}
{\cal D} = \sum_{p=0}^N {\cal R} d_p ,
\end{equation}
where $d_j = t_j\frac{\partial} {\partial t_j}$.
These  derivations extend to derivations of the Lie algebra
$\dot{\mathfrak{g}}\otimes {\cal R}$. Since,
$\hat\mathfrak{g}$ is the universal central extension of $\tilde\mathfrak{g}$,
we can lift these derivations to this universal central extension by  using a
result of
\cite{BM}. The action of vector fields on functions and
the Lie derivative action of vector fields on 1-forms leads to the following
action of ${\cal D}$ on $\hat\mathfrak{g}$:
\begin{equation}
\begin{aligned}
\label{action}
[t^{\bf m} d_j, g\otimes t^{\bf r}] &= r_j g\otimes t^{{\bf m}+{\bf r}} ,\\
[t^{\bf m} d_j, t^{\bf r} k_i] &= r_j t^{{\bf m}+{\bf r}}k_i + \delta_{ji}
\sum_{p=0}^N m_p t^{{\bf m}+{\bf r}}k_p .\end{aligned}
\end{equation}
 The formulas (\ref{action}) determine the Lie product on ${\cal D}$  up to a
${\cal K}$-valued 2-cocycle
$\tau \in H^2({\cal D}, {\cal K})$:
\begin{equation}\label{brackD}
[t^{\bf m} d_i , t^{\bf r} d_j] = r_i t^{{\bf m}+{\bf r}} d_j - m_j t^{{\bf
m}+{\bf r}} d_i + \tau( t^{\bf m} d_i , t^{\bf r} d_j).
\end{equation}
{}From the results of \cite{Dz}, any cocycle on ${\cal D}$ with values in
${\cal K}$ is a linear combination of
\begin{equation}\label{tau1}
\tau_1( t^{\bf m} d_i , t^{\bf r} d_j) = m_j r_i \sum_{p=0}^N r_p t^{{\bf
m}+{\bf r}} k_p
= - m_j r_i \sum_{p=0}^N m_p t^{{\bf m}+{\bf r}} k_p,
\end{equation}
and
\begin{equation}\label{tau2}
\tau_2( t^{\bf m} d_i , t^{\bf r} d_j) = m_i r_j \sum_{p=0}^N m_0 t^{{\bf
m}+{\bf r}} k_p.
\end{equation}
So we obtain the two-parametric family of algebras of \cite{BB}:
\[
\mathfrak{g}^{\cal D} = \mathfrak{g}_\tau^{\cal D}  = \hat\mathfrak{g}\oplus
{\cal D}=\dot{\mathfrak{g}} \otimes {\cal R} \oplus {\cal K} \oplus {\cal D},
\qquad
\text{where }\ \tau = \mu \tau_1 + \nu \tau_2.
\]
We denote by $\mathfrak{g}$ the following subalgebra of $\mathfrak{g}^{\cal
D}$:
\[
\mathfrak{g}= \hat\mathfrak{g}\oplus {\cal D}\oplus\bigoplus_{j=0}^N
\mathbb{C}d_j=\dot{\mathfrak{g}} \otimes {\cal R} \oplus {\cal K}
\oplus\bigoplus_{j=0}^N \mathbb{C}d_j .
\]

Note that the advantage of this larger Lie algebra is that
the center
of the algebra ${\mathfrak{g}} $ and ${\mathfrak{g}}^{\cal D} $is
finite-dimensional and is spanned
by $k_0, k_1, \ldots, k_N$, whereas the center of $\tilde{\mathfrak{g}} $ is
infinite dimensional.

Let $\Sigma=\{ \sigma_0,\sigma_1,\ldots,\sigma_N\}$ be a collection of finite
order automorphisms of $\dot{\mathfrak{g}}$. N.B., we do not assume that all
$\sigma_j$ are different and we allow $\sigma_j$ to be the identity. Let $n_j$
be the order of $\sigma_j $, i.e., $\sigma_j^{n_j}=1$ for the smallest
positive integer $n_j$ and denote by
$\epsilon_j=\exp \frac{2\pi i}{n_j}$. Then every $\sigma_j$ is diagonalizable
and one can decompose $\dot{\mathfrak{g}} $ in eigenspaces for  the eigenvalues
 $ \epsilon_j^k$, $k\in \mathbb{Z}/n_j\mathbb{Z}$. Assume from now on that all
$\sigma_j$, $0\le j\le N$ are simultaniously diagonalizable, i.e., one has the
following eigenspace decompostion of $\dot{\mathfrak{g}}$. Let ${\bf Z}$ be the
Cartesian product
\[
{\bf Z}=\mathbb{Z}/n_0\mathbb{Z}\times \mathbb{Z}/n_1\mathbb{Z}\times
\mathbb{Z}/n_2\mathbb{Z}\times\cdots\times \mathbb{Z}/n_N\mathbb{Z},
\]
then
\begin{equation}
\label{decomp}
\begin{aligned}
\dot{\mathfrak{g}} &=\bigoplus_{(k_0,k_1\cdots,k_N)\in{\bf Z}}
\dot{\mathfrak{g}}_{(k_0,k_1\cdots,k_N)}, \quad\text{where}\\[2mm]
\dot{\mathfrak{g}}_{(k_0,k_1\cdots,k_N)}&=\{ g\in
\dot{\mathfrak{g}}|\sigma_j(g)=\epsilon_j^{k_j} g\ \text{for all }0\le j\le
N\}.
\end{aligned}
\end{equation}

The Killing form $( \cdot , \cdot )$ is  ${\rm
Aut}\,\dot{\mathfrak{g}}$-invariant,
hence for every $0\le j\le N$ and all $x\in
\dot{\mathfrak{g}}_{(k_0,k_1\cdots,k_N)}$ and $y\in
\dot{\mathfrak{g}}_{(\ell_0,\ell_1\cdots,\ell_N)}$:
\[
(x,y)=(\sigma_j(x),(\sigma_j(y))=\epsilon_j^{k_j+\ell_j}(x,y).
\]
from which we conclude part (a) of the following Lemma:
\begin{lemma}
\label{L1}
(a) Let $( \cdot , \cdot )$ be the Killing form on $\dot{\mathfrak{g}}$, then
\[\left(\dot{\mathfrak{g}}_{(k_0,k_1\cdots,k_N)},
\dot{\mathfrak{g}}_{(\ell_0,\ell_1\cdots,\ell_N)}\right)=0\quad
\text{if }
(k_0+\ell_0,k_1+\ell_1\cdots,k_N+\ell_N)\ne (0,0,\cdots,0)\in{\bf Z}.
\]
(b) The subalgebra $\dot{\mathfrak{g}}_{(0,0,\cdots,0)}$ is reductive.
\end{lemma}
{\bf Proof} Part (b) of the Lemma is a direct consequence of part (a), the fact
that the Killing form is nondegenerate and Proposition 5 in \S 6.4 of
\cite{B}.\hfill{$\square$}
\

\

\noindent
This simple observation makes it possible to define twisted toroidal
subalgebras of a toroidal algebra. This construction, which we shall give now,
is similar to the one that produces the twisted affine Lie algebras  (see
\cite{Kac}, Chapter 8). But before we can do that, we will first introduce one
more notation. Let
${\bf m}= (m_0, m_1, ...,m_N) \in \mathbb{Z}^{N+1}$, then we write
\[
\overline{\bf m}=(\overline{m_0}, \overline{m_1}, ...,\overline{m_N})=
(m_0\ {\rm mod }n_0 , m_1\ {\rm mod }n_1, ...,m_N\ {\rm mod }n_n)\in{\bf Z}.
\]

Fix $\Sigma$, we define the subalgebra ${\tilde\mathfrak{g}}(\Sigma)$ of
${\tilde\mathfrak{g}}$ by
\begin{equation}
\label{twist1}
{\tilde\mathfrak{g}}(\Sigma)=\bigoplus_{{\bf m}\in\mathbb{Z}^{N+1}}
\dot{\mathfrak{g}}_{\overline{\bf m}}\otimes t^{\bf m}.
\end{equation}
Using Lemma \ref{L1}, one easily sees that one gets a subalgebra of
$\hat{\mathfrak{g}} $, which is a central extension of $\tilde{\mathfrak{g}} $,
if we add the subspace ${\cal K}(\Sigma)\subset {\cal K}$, which is spanned by
elements of the form
\[
t_0^{n_0m_0}t_1^{n_1m_1} \dots t_N^{m_Nm_N}k_i,\quad 0 \leq i \leq N,
\]
where of course the relation (\ref{exact}) still holds.
So define the following subalgebra of $\hat\mathfrak{g}$
\begin{equation}
\label{twist2}
\hat\mathfrak{g}(\Sigma)={\tilde\mathfrak{g}}(\Sigma)\oplus{\cal K}(\Sigma),
\end{equation}
with  Lie bracket on this algebra  still defined by (\ref{bracket}). We can
extend this twisted algebra with an algebra of derivations, however, except
when all automorphisms are the identity, not with ${\cal D}$, but with a
subalgebra of ${\cal D}$. Let
\begin{equation}
\label{twistD}
{\cal D}(\Sigma) = \sum_{p=0}^N {\cal R}(\Sigma) d_p ,\quad\text{where }\
{\cal R}(\Sigma)=\mathbb{C}[t_0^{\pm n_0},t_1^{\pm n_1}
,\ldots t_N^{\pm n_N}],
\end{equation}
define a subalgebra $\mathfrak{g}^{\cal D}(\Sigma)$ of $\mathfrak{g}^{\cal D}$
and a subalgebra $\mathfrak{g}(\Sigma)$ of $\mathfrak{g}$ by
\[
\begin{aligned}
\mathfrak{g}^{\cal D}(\Sigma) &= \mathfrak{g}_\tau^{\cal D} (\Sigma)=
\hat\mathfrak{g}(\Sigma)\oplus {\cal
D}(\Sigma)=\tilde\mathfrak{g}(\Sigma)\oplus{\cal K}(\Sigma) \oplus{\cal
D}(\Sigma),\\
\mathfrak{g}(\Sigma)&= \hat\mathfrak{g}(\Sigma) \oplus\bigoplus_{j=0}^N
\mathbb{C}d_j=\tilde\mathfrak{g}(\Sigma)\oplus{\cal K}(\Sigma)
\oplus\bigoplus_{j=0}^N \mathbb{C}d_j,
\end{aligned}
\]
where the Lie bracket is still defined by (\ref{bracket}), (\ref{action}) and
(\ref{brackD}).

Let $\dot\mathfrak{h}_{\overline{\bf 0}}$ be the Cartan subalgebra of
$\dot\mathfrak{g}_{(\overline 0,\overline 0,\ldots, \overline 0)}$,
then
\[
\mathfrak{h}=\dot\mathfrak{h}_{\overline{\bf 0}}\oplus \mathbb{C}d_0
\oplus \mathbb{C}d_1 \oplus\cdots\oplus  \mathbb{C}d_N \oplus \mathbb{C}k_0
\oplus \mathbb{C}k_1 \oplus\cdots\oplus \mathbb{C}k_N
\]
is the Cartan subalgebra of $\mathfrak{g}^{\cal D}(\Sigma)$ and
$\mathfrak{g}(\Sigma)$. We extend $\lambda\in
\dot\mathfrak{h}_{\overline{\bf 0}}^*$ to a linear  function on $\mathfrak{h}$
by setting $\lambda(d_i)=\lambda(k_i)=0$, for all $0\le i\le N$.
Denote by
$\delta_i, \kappa_i$ the linear function on $\mathfrak{h}$ defined by
\[
\begin{aligned}
\delta_i(\dot\mathfrak{h}_{\overline{\bf 0}})=0,&\quad
\delta_i(d_j)=\delta_{ij}, \quad \delta_i(\kappa_j)=0,\\
\kappa_i(\dot\mathfrak{h}_{\overline{\bf 0}})=0,&\quad
\kappa_i(d_j)=0, \quad \kappa_i(\kappa_j)=\delta_{ij}.
\end{aligned}
\]
The Killing form restricted to $\dot\mathfrak{h}_{\overline{\bf 0}}$ remains
nondegenerate and  can be extended to a nondegenerate symmetric bilinear form
on $\mathfrak{h}$
\[
(k_i,k_j)=0,\quad (d_i,k_j)=\delta_{ij},\quad (d_i,d_j)=0,\quad
(k_i,\dot\mathfrak{h}_{\overline{\bf 0}})=(d_i,\dot\mathfrak{h}_{\overline{\bf
0}})=0.
\]
This form defines an isomorphism $\mathfrak{h}\to\mathfrak{h}^*$ by
\[
\nu(h)(h')=(h,h'),\qquad h,h'\in\mathfrak{h}_{\overline{\bf 0}}
\]
and hence a bilinear form on $\mathfrak{h}^*$, viz.,
\[
(\alpha,\beta)=(\nu^{-1}(\alpha),\nu^{-1}(\beta)).
\]
One thus has
\[
(\kappa_i,\kappa_j)=0,\quad (\delta _i, \kappa_j)=\delta_{ij},\quad
(\delta_i,\delta_j)=0,\quad
(\kappa_i,\dot\mathfrak{h}_{\overline{\bf
0}}^*)=(\delta_i,\dot\mathfrak{h}_{\overline{\bf 0}}^*)=0.
\]
Then, $\mathfrak{g}(\Sigma)$ decomposes with respect to $\mathfrak{h}^*$ into
\[
\mathfrak{g}(\Sigma)=
\bigoplus_{\alpha\in\mathfrak{h}^*}\mathfrak{g}(\Sigma)_\alpha.
\]
Let
\[
\Delta=\{ \alpha \in \mathfrak{h}^*|
\mathfrak{g}(\Sigma)_\alpha\ne \{ 0\}\},
\]
be the set of roots of $\mathfrak{g}(\Sigma)$
then we have the following root space decomposition
\[
\mathfrak{g}(\Sigma)=\mathfrak{h}\oplus\bigoplus_{\alpha\in
\Delta}\mathfrak{g}(\Sigma)_\alpha.
\]
The connection with the extended affine root system of Section \ref{EARS} is as
follows.
The linear space $V$ in Definition \ref{D1}, is the subspace
\[
V= \dot\mathfrak{h}_{\overline{\bf 0}}^*\oplus\bigoplus_{i=0}^N
\mathbb{C}\delta_i,
\]
and the bilinear form of the definition is the restriction  $(\cdot,\cdot)_V$
of $(\cdot,\cdot)$ to $V$. One can decompose $\mathfrak{g}(\Sigma)$  with
respect to $V$ into
\[
\mathfrak{g}(\Sigma)=\bigoplus_{\alpha\in V}\mathfrak{g}(\Sigma)_\alpha.
\]
Let
\[
\overline R=\{ \alpha \in V| \mathfrak{g}(\Sigma)_\alpha\ne \{ 0\}\},
\]
then
\[
\overline R=R\cup R^0,\qquad \text{where}\quad
R=\{\alpha\in \overline R | (\alpha,\alpha)_V\ne 0\}
\quad\text{and}\quad
R^0=\overline R\cap V^0
\]
and thus
\[
\mathfrak{g}(\Sigma)=\mathfrak{h}\oplus\bigoplus_{\alpha\in
R}\mathfrak{g}(\Sigma)_\alpha\oplus\bigoplus_{\alpha\in
R^0}\mathfrak{g}(\Sigma)_\alpha.
\]

Note that there is one problem, for general $\Sigma$, it is not clear that the
set $R$ satisfies the axioms of Definition \ref{D1}.
In the next section we choose $N=1$ and list pairs $\dot\mathfrak{g}$, $\Sigma$
which give the extended affine Lie algebras that correspond to the extended
affine root systems of Saito \cite{S1}, i.e., to the ones that were presented
in section \ref{EARS}.

\section{Bi-affine Lie Algebras}
\label{BLA}

In this section we construct the bi-affine Lie algebras, i.e., the twisted
2-toroidal Lie algebras corresponding to the the extended affine root systems
of Saito \cite{S1}, which were presented in Section \ref{EARS}. So we assume
from now on in this section that $N=1$. In most cases we will explain how we
realize
$\dot\mathfrak{g}$, this will however not always be the same, e.g. the Lie
algebra of type $D_\ell$ will be realized in different ways.

\subsection{Type $X_\ell^{(1,1)}$}
Bi-affine Lie algebras of type $X_\ell^{(1,1)}$, can be easily constructed.
One takes in the construction of section \ref{TA} for  $\dot\mathfrak{g}$ the
simple Lie algebra of type $X_\ell$ and as automorphisms
$\sigma_0=\sigma_1=\text{id}$.

\subsection{Type $X_\ell^{(1,t)}$}
\label{1t}
The description of the bi-affine Lie algebras of type $X_\ell^{(1,t)}$ is also
easy. For $X_\ell$ equal to $B_\ell$, $C_\ell$, $F_4$ and $G_2$, one takes as
$\dot\mathfrak{g}$ the simple Lie algebra of type  $D_{\ell+1}$, $A_{2\ell-1}$,
$E_6$, $D_4$, respectively. One chooses for $\sigma_0$ the identity and for
$\sigma_1$ the automorphisms, described in \S 8 of \cite{Kac}, which are
induced by a diagram automorphism. The order of $\sigma_2$ is $t$, which is
equal to 2, 2, 2, 3, respectively.

\subsection{Type $X_\ell^{(t,t)}$}
\label{T3}
We start this subsection with the $X_\ell=B_\ell$ and $t=2$. We take as
$\dot\mathfrak{g}$ the Lie algebra of type $D_{\ell+2}$. Let $M_n$ be the
linear space of all complex $n\times n$-matrices. We realize $\dot\mathfrak{g}$
as
\begin{equation}
\label{so}
\dot\mathfrak{g}=\{ X\in M_{2\ell+4} | X^T=-X \},
\end{equation}
where $X^T$ stands for the transposed of the matrix $X$.
Let $E_{ij}$ be the matrix with a 1 on the $(i,j)$-th entry and zeros
elsewhere.
Now choose $\Sigma$ as follows
\[
\begin{aligned}
\sigma_0=&{\rm Ad}
\left(-E_{2\ell+3,2\ell+3}-E_{\ell+4,\ell+4}+\sum_{i=1}^{2\ell+2} E_{ii}\right)
,\\
\sigma_1=&{\rm Ad}
\left(-E_{2\ell+2,2\ell+2}+E_{2\ell+3,2\ell+3}
-E_{2\ell+4,2\ell+4}+\sum_{i=1}^{2\ell+1} E_{ii}\right).
\end{aligned}
\]
The subalgebra  $\dot\mathfrak{g}_{(\overline 0,\overline 0)}$ is the simple
Lie algebra of type $B_\ell$. All three other spaces consist of the direct sum
of a 1-dimensional trivial and a $2\ell+1$-dimensional irreducible
representation of $B_\ell$.

Next we take $X_\ell=C_\ell$ and $t=2$. In this case $\dot\mathfrak{g}$ is the
Lie algebra of type $D_{2\ell}$, which we realize as
\begin{equation}
\label{typeD}
\dot\mathfrak{g}=\left \{
\left(\begin{array}{cc} a&b\\ c&-a^T
\end{array}\right)\in M_{4\ell}| b^T=-b,\ c^T=-c
\right\}.
\end{equation}
\label{s12}
The automorphisms are defined as follows
\begin{equation}
\begin{aligned}
\sigma_0=&{\rm Ad} \left(\sum_{i=1}^{4\ell}(-)^i E_{i,4\ell+1-i}\right)
,\\
\sigma_1=&{\rm Ad} \left(\sum_{i=1}^{2\ell}(-)^i E_{i,4\ell+1-i}-(-)^i
E_{2\ell+i,2\ell+1-i}\right).
\end{aligned}
\end{equation}
Here $\dot\mathfrak{g}_{(\overline 0,\overline 0)}$ is the simple Lie algebra
of type $C_\ell$. All three other spaces consist of the direct sum of a
1-dimensional trivial and a $2\ell^2-\ell-1$-dimensional irreducible
representation of $B_\ell$.

Assum now that $X_\ell=F_4$, then $t=2$. Now $\dot\mathfrak{g}$ is the Lie
algebra of type $E_7$, and assume that the roots are labeled as "Planche VI" in
$\cite{B2}$. Let $e_i$, $f_i$, $1\le i\le 7$ be the Chevalley generators
corresponding to these roots. We define both automorphisms on these generators.
\[
\begin{aligned}
\sigma_0(e_i)=e_i,&\quad \sigma_0(f_i)=f_i\quad  \text{for } i=2,4,\\
\sigma_0(e_1)=e_6,&\quad \sigma_0(f_1)=f_6,\\
\sigma_0(e_3)=e_5,&\quad \sigma_0(f_3)=f_5,\\
\sigma_0(e_5)=e_3,&\quad \sigma_0(f_5)=f_3,\\
\sigma_0(e_6)=e_1,&\quad \sigma_0(f_6)=f_1,\\
\sigma_0(e_7)=e_{-\theta}&
=[f_1f_3f_4f_2f_5f_4f_6f_5f_3f_4f_2f_1f_3f_4f_5f_6f_7],
\quad \sigma_0(f_7)=e_\theta,
\end{aligned}
\]
where $[ab\ldots cd]$ stands for $[a,[b.\ldots,[c,d]]]\ldots]$. Here $e_\theta$
is a root vector corresponding to the highest root and $e_{-\theta}$ is a root
vector corresponding to the lowest root. Of course $e_\theta$ has to be
normalized in such a way that $\sigma_0$ is an automorphism of order 2.
\[
\begin{aligned}
\sigma_1(e_i)=e_i,&\quad \sigma_1(f_i)=f_i\quad  \text{for } i\ne 7\\
\sigma_1(e_7)=-e_7,&\quad \sigma_1(f_7)=-f_7.
\end{aligned}
\]
Then with respect to these two automorphisms $\dot\mathfrak{g}$ splits into the
following spaces
$\dot\mathfrak{g}_{(\overline 0,\overline 0)}$ is the simple Lie algebra of
type $F_4$ and all the spaces $\dot\mathfrak{g}_{(\overline 0,\overline 1)}$,
$\dot\mathfrak{g}_{(\overline 1,\overline 0)}$ and
$\dot\mathfrak{g}_{(\overline 1,\overline 1)}$ is the direct sum of a
1-dimensional trivial representation  and the 26-dimensional irreducible
representation of $F_4$. This was checked using
\emph{ {GAP -- Groups, Algorithms, and Programming}} \cite{GAP4}.

Finally, we describe $G_2^{(3,3)}$. For this case $\dot\mathfrak{g}$ is the Lie
algebra of type $E_6$. We define the automorphisms again  on the Chevalley
generators,where the roots are numbered as in "Planche V" of $\cite{B2}$.
Let $\theta$ be the highest root of $E_6$, now define
\[
\begin{aligned}
\sigma_0(e_1)=e_6,&\quad \sigma_0(f_1)=f_6,\\
\sigma_0(e_2)=e_3,&\quad \sigma_0(f_2)=f_3,\\
\sigma_0(e_3)=e_5,&\quad \sigma_0(f_3)=f_5,\\
\sigma_0(e_4)=e_4,&\quad \sigma_0(f_4)=f_4,\\
\sigma_0(e_5)=e_2,&\quad \sigma_0(f_5)=f_2,\\
\sigma_0(e_6)=e_{-\theta}&=[f_2f_4f_5f_3f_4f_2f_6f_5f_4f_3f_1],\quad
\sigma_0(f_7)=e_\theta,
\end{aligned}
\]
which is an automorphism of order 3. Let $\omega=e^{\frac{2\pi i}{3}}$, define
$\sigma_1$ as follows:
\[
\begin{aligned}
\sigma_1(e_i)=e_i,&\quad \sigma_1(f_i)=f_i\quad  \text{for } i\ne 1,6,\\
\sigma_1(e_i)=\omega e_i,&\quad \sigma_1(f_i)=\omega^2f_i.\quad  \text{for } i=
1,6.
\end{aligned}
\]
The 78-dimensional Lie algebra $E_6$ splits with respect to these automorphisms
as follows: $\dot\mathfrak{g}_{(\overline 0,\overline 0)}$ is the simple Lie
algebra of type $G_2$, all 8 other spaces $\dot\mathfrak{g}_{(\overline
i,\overline j)}$, $0\le i,j\le 2$, $(i,j)\ne (0,0)$ is the direct sum of the
irreducible 7-dimensional representation and a 1-dimensional trivial
representation of $G_2$.
This example was also checked using \emph{ {GAP}} \cite{GAP4}.

\subsection{Type $A_1^{(1,1)*}$}

We will not describe the bi-affine Lie algebra of this type here, but instead
we will obtain it as a special case of the construction of the next section.

\subsection{Type $B_\ell^{(2,2)*}$}

The construction of the bi-affine Lie algebra of type $B_\ell^{(2,2)*}$ goes as
follows. The algebra $\dot\mathfrak{g}$ is the simple Lie algebra of type
$B_{\ell+1}$, which we realize as the space of anti-symmetric matrices (cf.
(\ref{so})):
\[
\dot\mathfrak{g}=\{ X\in M_{2\ell+3} | X^T=-X \}.
\]
The automorphisms are taken as follows:
\[
\begin{aligned}
\sigma_0=&{\rm Ad} \left(-E_{2\ell+3,2\ell+3}+\sum_{i=1}^{2\ell+2}
E_{ii}\right)
,\\
\sigma_1=&{\rm Ad}
\left(-E_{2\ell+2,2\ell+2}+E_{2\ell+3,2\ell+3}+\sum_{i=1}^{2\ell+1}
E_{ii}\right).
\end{aligned}
\]
If we take $\ell=1$, we obtain the bi-affine Lie algebra of type
$A_1^{(1,1)*}$.

Here $\dot\mathfrak{g}_{(\overline 0,\overline 0)}$ is the simple Lie algebra
of type $B_\ell$. The spaces $\dot\mathfrak{g}_{(\overline 0,\overline 1)}$ and
$\dot\mathfrak{g}_{(\overline 1,\overline 0)}$
consist of   a $2\ell+1$-dimensional irreducible representation of $B_\ell$
and $\dot\mathfrak{g}_{(\overline 1,\overline 1)}$ is one dimensional.

\subsection{Type $C_\ell^{(1,1)*}$}

To obtain the  bi-affine Lie algebra of type $C_\ell^{(1,1)*}$ we choose
$\dot\mathfrak{g}$ to be the simple Lie algebra of type $C_{2\ell}$.
This Lie algebra is realized as follows (cf. (\ref{typeD})):
\[\dot\mathfrak{g}=\left \{
\left(\begin{array}{cc} a&b\\ c&-a^T
\end{array}\right)\in M_{4\ell}| b^T=b,\ c^T=c
\right\}.
\]
The automorphisms $\sigma_0$ and $\sigma_1$ are defined by (\ref{s12}).
Now, $\dot\mathfrak{g}_{(\overline 0,\overline 0)}$ is the simple Lie algebra
of type $C_\ell$.Both spaces $\dot\mathfrak{g}_{(\overline 0,\overline 1)}$ and
$\dot\mathfrak{g}_{(\overline 1,\overline 0)}$ are the adjoint representation
and
$\dot\mathfrak{g}_{(\overline 1,\overline 1)}$
consist of  the direct sum of a $2\ell^2-\ell-1$-dimensional irreducible
representation and a one dimensional one.

\subsection{Type $BC_\ell^{(2,1)}$}

This type can be described in the same way as the examples of Section \ref{1t}.
Here $\dot\mathfrak{g}$ is the Lie algebra of type $A_{2\ell}$. Let $\sigma_0$
be the automorphism, described  in \S 8 of \cite{Kac}, which produces the
affine Lie algebra of type $A_{2\ell}^{(2)}$, for $\sigma_1$ we choose the
identity.

\subsection{Type $BC_\ell^{(2,2)}(1)$}

The bi-affine algebra of type  $BC_\ell^{(2,1)}$ is constructed by taking as
$\dot\mathfrak{g}$ the simple Lie algebra of type $D_{2\ell+1}$, realized as
\[
\dot\mathfrak{g}=\{ X\in M_{4\ell+2} | X^T=-X \}.
\]
In this case $\Sigma $ is given by
\[
\begin{aligned}
\sigma_0=&{\rm Ad} \left(\sum_{i=1}^{4\ell+2}(-)^i E_{i,4\ell+3-i}\right)
,\\
\sigma_1=&{\rm Ad} \left(\sum_{i=1}^{2\ell+1}
E_{i,i}-E_{2\ell+1+i,2\ell+1+i}\right).
\end{aligned}
\]
In this example is  $\dot\mathfrak{g}_{(\overline 0,\overline 0)}$  the simple
Lie algebra of type $B_\ell$. Both spaces $\dot\mathfrak{g}_{(\overline
0,\overline 1)}$ and
$\dot\mathfrak{g}_{(\overline 1,\overline 1)}$ are the adjoint representation
and
$\dot\mathfrak{g}_{(\overline 1,\overline 0)}$
consist of  the direct sum of a $2\ell^2+3\ell$-dimensional irreducible
representation and a one dimensional one.

\subsection{Type $BC_\ell^{(2,2)}(2)$}
\label{T9}

Type $BC_\ell^{(2,2)}(2)$ can be obtained by choosing $\dot\mathfrak{g}$ the
simple Lie algebra of type $A_{2\ell+1}$. We realize this Lie algebra in the
usual way as the traceless $(2\ell+2)\times(2\ell+2)$-matrices. $\sigma_0$ is
the Cartan involution
\[
\begin{aligned}
\sigma_0(X)&=-X^T\quad  \text{ and} \\
\sigma_1&={\rm Ad} \left(-E_{2\ell+2,2\ell+2}+\sum_{i=1}^{2\ell+1}
E_{i,i}\right).
\end{aligned}
\]
For the decomposition of this Lie algebra we refer the reader to the second
example of Section \ref{QLA}.

\subsection{Type $BC_\ell^{(2,4)}$}

Finally, the  bi-affine Lie algebra of type $BC_\ell^{(2,4)}$ is constructed as
follows. The algebra $\dot\mathfrak{g}$  is the Lie algebra of type
$D_{2\ell+2}$,
which we realize as in (\ref{typeD}), but then with $4\ell$ replaced by
$4\ell+4$. The automorphism $\sigma_0$ is the involution
\[
\sigma_0={\rm Ad}
\left(E_{2\ell+2,4\ell+4}-E_{4\ell+4,2\ell+2}+\sum_{j=1}^{2\ell+1}
(-)^jE_{4\ell+4-j,j}-(-)^jE_{2\ell+2-j,2\ell+2+j}\right).
\]
The other automorphism is an automorphism of order 4:
\[
\sigma_1={\rm Ad}
\left(iE_{2\ell+2,2\ell+4}+iE_{2\ell+4,2\ell+2}+\sum_{j=1}^{2\ell+1}
E_{j,j}-E_{2\ell+2+j,2\ell+2+j}\right).
\]
The decomposition of $D_{2\ell+2}$ with respect to the automorphisms is the
most complicated one. The subalgebra $\dot\mathfrak{g}_{(\overline 0,\overline
0)}$  the simple Lie algebra of type $B_\ell$. The spaces
$\dot\mathfrak{g}_{(\overline i,\overline j)}$, where $(i,j)= (0,1),\ (1,1),\
(0,3)$ or $(1,3)$ are all irreducible $2\ell+1$ dimensional representations of
$B\ell$. Next, $\dot\mathfrak{g}_{(\overline 1,\overline 0)}$ is the direct sum
of a one dimensional  and the $2\ell^2+3\ell$-dimensional irreducible
representation of $B_\ell$;
$\dot\mathfrak{g}_{(\overline 0,\overline 2)}$ is the ad-module (the adjoint
representation) and
$\dot\mathfrak{g}_{(\overline 1,\overline 2)}$ is again the direct sum of two
irreducible modules, viz.  the adjoint representation and a one dimensional
trivial one.

This produces all bi-affine Lie algebras related to Saito's list. The ones of
type $X_\ell^{(t,1)}$ are isomorphic to the ones of type $X_\ell^{(1,t)}$.
The former can be constructed by interchanging the automorphisms $\sigma_1$ and
$\sigma_2$ in the construction of latter. Except for two cases, viz.,
$A_1^{(1,1)*}$ and $C_\ell^{(1,1)*}$, the tier numbers, introduced by Saito,
which are the upper indeces in $X_\ell^{(s,t)}$, equal the order of the
automorphims used to construct the bi-affine lie algebras.

\section{Some quasi-simple Lie algebras}
\label{QLA}

In 1990  Hoegh-Krohn and Torresani \cite{HT} classified  and constructed
certain quasi-simple Lie algebras. These are characterized by the existence of
a finite-dimensional Cartan subalgebra, a nondegenerate invariant symmetric
bilinear form and nilpotent root spaces attached to non-isotropic roots. They
derive a classification for the possible irreducible elliptic quasi-simple root
systems. Obviously they were not aware of the existence of the paper of Saito
\cite{S1}, which appeared 5 years earlier. In the case when the nullity is
equal to 2, their list lacked some of the cases Saito obtained. According to
the introduction of \cite{AA}, this was caused by the fact that they assumed or
erroneously concluded that the theory governing the isotropic roots was based
on lattices. Two of their examples were the root systems
\begin{equation}
\label{Rex}
\begin{aligned}
R(X_\ell)^{(t,t,\ldots,t,1,1,\ldots ,1)}=&\{\alpha+\sum_{j=0}^n k_j\delta_j
|\alpha \in R(X_\ell)\ \text{short},\ k_j\in\mathbb{Z}\}\\
\ &\cup \{\alpha+\sum_{j=0}^m tk_j\delta_j +\sum_{j=m+1}^n k_j\delta_j |\alpha
\in R(X_\ell)\ \text{long},\ k_j\in\mathbb{Z}\},
\end{aligned}
\end{equation}
with
$X_\ell=B_\ell,\ C_\ell$ and $F_4$, and  $X_\ell=G_2$ and $t= 2$, 2, 2, 3,
respectively and $0\le m\le n$;
\begin{equation}
\label{Rex2}
\begin{aligned}
R&(BC_\ell)^{(2,2,\ldots,2,1,1,\ldots ,1)}(2)=\{\alpha+\sum_{j=0}^n k_j\delta_j
|\alpha \in R(B_\ell)\ \text{short},\ k_j\in\mathbb{Z}\}\\
\ &\cup\{\alpha+k_0\delta_0+\sum_{j=1}^{m+1} 2k_j\delta_j +\sum_{j=m+2}^{n}
k_j\delta_j|\alpha \in R(B_\ell)\ \text{long},\ k_j\in\mathbb{Z}\}\\
\ &\cup \{\alpha+(2k_0+1)\delta_0+\sum_{j=1}^{m+1}2k_j\delta_j
+\sum_{j=m+2}^{n} k_j\delta_j|\alpha \in R(C_\ell)\ \text{long},\
k_j\in\mathbb{Z}\},
\end{aligned}
\end{equation}
with $0\le m< n$.

Using the idea's of Sections \ref{1t}, \ref{T3}, \ref{T9}
we construct the corresponding extended affine algebra related to (\ref{Rex})
for $X_\ell=B_\ell$ and to (\ref{Rex2}). This shows that the approach of this
paper, not only works for $N=1$, but that it at least also produces some
families of $N$-affine Lie algebras.

For the example related to (\ref{Rex}) we take as  $\dot\mathfrak{g}$ the Lie
algebra of type $D_{\ell+2^m}$.
We again realize this Lie algebra as the complex space of anti-symmetric
$(2\ell+2^{m+1})\times(2\ell+2^{m+1})$-matrices (cf. (\ref{so})). Let
$J_0,J_1,\ldots
J_m$ be the following matrices
\[
\begin{aligned}
J_0&=\sum_{i=1}^{2\ell}E_{ii}+\sum_{j=1}^{2^m}
\left(E_{2\ell+2j-1,2\ell+2j-1}-E_{2\ell+2j,2\ell+2j}\right),\\[2mm]
J_1&=\sum_{i=1}^{2\ell}E_{ii}+\sum_{j=1}^{2^{m-1}}
\left(E_{2\ell+4j-3,2\ell+4j-3}+
E_{2\ell+4j-2,2\ell+4j-2}-E_{2\ell+4j-1,2\ell+4j-1}
+E_{2\ell+4j,2\ell+4j}\right),\\[2mm]
J_2&=\sum_{i=1}^{2\ell}E_{ii}+\sum_{j=1}^{2^{m-2}}\left(
E_{2\ell+8j-7,2\ell+8j-7}+E_{2\ell+8j-6,2\ell+8j-6}+E_{2\ell+8j-5,2\ell+8j-5}
+E_{2\ell+8j-4,2\ell+8j-4}
\right.\\[1mm]
\ &\qquad\qquad\qquad \left.
-E_{2\ell+8j-3,2\ell+8j-3}-E_{2\ell+8j-2,2\ell+8j-2}
-E_{2\ell+8j-1,2\ell+8j-1}-E_{2\ell+8j,2\ell+8j}\right),\\
\ &\ \vdots\\
J_m&=\sum_{i=1}^{2\ell+2^m}E_{ii}-\sum_{j=1}^{2^{m}}
E_{2\ell+2^m+j,2\ell+2^m+j},
\end{aligned}
\]
then we define $\Sigma=(\sigma_0,\sigma_1,\ldots\sigma_n)$ by
\[
\sigma_k=
\begin{cases}{\rm Ad}(J_k)&\quad \text{for }0\le k\le m,\\
{\rm id}&\quad \text{for }k>m.
\end{cases}
\]
Then $\dot\mathfrak{g}_{(\overline 0,\overline 0,\ldots,\overline 0)}$ is the
simple Lie algebra of type $B_\ell$ and all the $2^{m+1}-1$ other spaces
$\dot\mathfrak{g}_{(\overline i_0,\overline i_1\ldots \overline i_m,\overline
0,\ldots,\overline 0)}$ consist of the direct sum of the $2\ell+1$-dimensional
irreducible representation together with $2^m-1$ trivial 1-dimensional ones.
Clearly the Lie algebra $\mathfrak{g}(\Sigma)$ corresponds to the root system
for
$X_\ell=B_\ell$
which is given in (\ref{Rex}).

For the second example, the one related to (\ref{Rex2}), we take as
$\dot\mathfrak{g}$ the simple Lie algebra of type $A_{2\ell+2^{m+1}-1}$, which
can be realized in the usual way as complex traceless $2\ell+2^{m+1}\times
2\ell+2^{m+1}$-matrices. For $\sigma_0$ we again take as \S \ref{T9} the Cartan
involution, i.e., $\sigma_0(X)=-X^T$. All the other automorphisms are defined
as follows:
\[
\sigma_k=
\begin{cases}{\rm Ad}(J_{k-1})&\quad \text{for }1\le k\le m+1,\\
{\rm id}&\quad \text{for }k>m+1.
\end{cases}
\]
The algebra  $\dot\mathfrak{g}_{(\overline 0,\overline 0,\ldots,\overline 0)}$
is the simple Lie algebra of type $B_\ell$. The space
$\dot\mathfrak{g}_{(\overline 1,\overline 0,\ldots,\overline 0)}$ consists of
the direct sum of a $2\ell^2+3\ell$-dimensional irreducible representation of
$B_\ell$ together with $2^m-1$ 1-dimensional trivial representations.
All other $2^{m+2}-2$ spaces $\dot\mathfrak{g}_{(\overline i_0,\overline
i_1,\ldots,\overline i_{m+1},\overline 0,\ldots,\overline 0)}$ cosist of the
direct sum of the $2\ell+1$ irreducible representation of $B_\ell$
and $2^m-1$ trivial 1-dimensional representations. This leads to a Lie algebra
 $\mathfrak{g}(\Sigma)$, whose  root system restricted to $V$ is given by
(\ref{Rex2}).

\end{document}